\documentclass [11pt,oneside,a4paper,mathscr]{amsart}
\usepackage {amsmath,amssymb,euscript}

\newtheorem{thm}{Theorem}[section]
\newtheorem{lemma}[thm]{Lemma}
\newtheorem{cor}[thm]{Corollary}
\newtheorem{prop}[thm]{Proposition}

\theoremstyle{definition}
\newtheorem{dfn}[thm]{Definition}
\newtheorem{remark}[thm]{Remark}
\newtheorem*{ack}{Acknowledgements}
\newtheorem*{nota}{Notation}

\theoremstyle{remark}
\newtheorem*{pf}{Proof}

\numberwithin{equation}{subsection}

\let\tensor=\otimes

\newcommand{\KK}{\mathcal K}
\newcommand{\Ltensor}{\mathbin{\buildrel{\mathbf L}\over{\tensor}}}
\newcommand{\C}{{\mathbb C}}
\newcommand{\I}{\mathcal I}
\newcommand{\F}{{\mathcal F}}
\newcommand{\Z}{{\mathbb Z}}
\newcommand{\PP}{{\mathcal P}}
\newcommand{\OO}{{\mathscr O}}
\newcommand{\id}{\operatorname{Id}}
\newcommand{\Hilb}{\operatorname{Hilb}}
\newcommand{\Pic}{\operatorname{Pic}} 
\newcommand{\Ext}{\operatorname{Ext}}

\newcommand{\Hom}{\operatorname{Hom}}
\newcommand{\ch}{\operatorname{ch}}
\newcommand{\Spec}{\operatorname{Spec}}
\newcommand{\D}{{D}}
\renewcommand{\leq}{\leqslant}
\renewcommand{\geq}{\geqslant}

\newcommand{\bda}{\big\downarrow}
\newcommand{\lra}{\longrightarrow}
\newcommand{\lRa}[1]{\>{\buildrel {#1}\over\longrightarrow}\>}
\newcommand{\M}{{\mathcal M}}
\newcommand{\NS}{\operatorname{NS}}
\newcommand{\dual}{\vee}
\newcommand{\EE}{\mathcal E}
\newcommand{\FF}{\mathcal F}
\newcommand{\QQ}{\mathcal Q}
\newcommand{\N}{\mathbb N}
\newcommand{\R}{\mathbf R}
\renewcommand{\L}{\mathbf L}
\newcommand{\ellsur}{X\lRa{\pi} C}
\newcommand{\isom}{\cong}
\newcommand{\J}{J}
\newcommand{\embed}{\hookrightarrow}
\newcommand{\HH}{\mathscr H}
\newcommand{\lHom}{\operatorname{{\mathcal H}om}}
\newcommand{\td}{\operatorname{td}}
\newcommand{\col}[2]{\left( \begin{array}{c} #1 \\ #2 \end{array} \right)}
\newcommand{\mat}[4]{\left( \begin{array}{cc} #1 & #2 \\ #3 & #4
\end{array} \right)}
\newcommand{\m}{\lambda}
\newcommand{\De}{\Lambda}
\newcommand{\NN}{\mathscr N}
\newcommand{\U}{\mathscr U}
\newcommand{\V}{\mathscr V}
\newcommand{\W}{\mathscr W}

\newcommand{\oo}{\circ}

\begin{document}

\title[]{Fourier-Mukai transforms for elliptic surfaces}
\author[]{Tom Bridgeland}
\date{\today}
\subjclass{14D20, 14J27, 14J60}

\begin{abstract}
We compute a large number of moduli spaces of stable bundles on a
general algebraic elliptic surface using a new class of relative Fourier-Mukai transforms.
\end{abstract}

\maketitle

\section{Introduction}

In recent years moduli
spaces of
stable bundles on projective surfaces have been extensively
studied. Although various
general results are known, few of these spaces have been described explicitly. One important development was
S. Mukai's discovery [11] of a
transform which allowed him to compute some of the moduli spaces on
Abelian varieties [14]. This technique has recently been extended to cover K3
surfaces [4], [13], [17]. Here we introduce similar transforms for elliptic surfaces.

Much is already known about moduli of rank 2 bundles on elliptic
surfaces thanks to work by R. Friedman [7], [8], amongst others.
This research has led
to several important results, including the smooth classification of
elliptic surfaces. In this paper we use Mukai's techniques to
study bundles of higher rank.

\subsection{}
Let $\ellsur$ be a relatively minimal algebraic elliptic surface over $\C$.
Given a sheaf
$E$ on $X$ we write its Chern class as a triple
$$(r(E),c_1(E),c_2(E))\in\Z\times\NS(X)\times\Z.$$
Here $\NS(X)$ is the
N\'eron-Severi group of $X$, i.e. the subgroup of $H^2(X,\Z)$ generated
by the Chern classes of line
bundles on $X$. We denote the element of $\NS(X)$ corresponding to a fibre of $\pi$ by $f$.

Let $(r,\De,k)$ be a triple as above, and assume that $r>1$ is coprime to
$\De\cdot f$. Then, as Friedman observed, there exist polarizations of $X$
with respect to which a torsion-free sheaf $E$ of Chern class
$(r,\De,k)$ is stable iff the
restriction of $E$ to the general fibre of $\pi$ is stable. For these
polarizations, and sheaves of this Chern class,
the notions of Gieseker stability, $\mu$-semi-stability
and $\mu$-stability all coincide. Taking such
a polarization, we define $\M=\M_X(r,\De,k)$ to be the moduli space of stable
torsion-free sheaves on $X$ of Chern class $(r,\De,k)$.

An argument of Friedman's shows that the projective scheme $\M$
is smooth of
dimension $\dim(\Pic^{\oo}(X))+2t$, where
\begin{equation}
\label{first}
2t=2rk-(r-1)\De^2-(r^2-1)\chi(\OO_X).
\end{equation}
If $t<0$ then $\M$ is empty, so we shall assume that $t$ is
non-negative. Our main result (Theorem 1.1 below) states that $\M$ is irreducible and
birationally equivalent to $\Pic^{\oo}(Y)\times\Hilb^t(Y)$, where $Y$
is another elliptic surface over $C$.

\subsection{}
Let us define $\m_X$ to be the
highest common factor of the fibre degrees of sheaves on $X$. Equivalently, $\m_X$ is the smallest positive
integer such that $\ellsur$ has a holomorphic $\m_X$-multisection.

Given a pair of
integers $a>0$ and $b$ such that $a\m_X$ is coprime to $b$, we define
an elliptic surface $\J_X(a,b)$ over $C$, whose fibre over a point
$p\in C$ is canonically identified with (a component of) the moduli
space of rank $a$, degree $b$, stable sheaves on the fibre $X_p$.

To do this, take a polarization of $X$ of fibre degree
coprime to $b$, and let $\M(X/C)\to C$ denote the relative moduli
space of stable pure-dimension 1 sheaves on the fibres of $\pi$
(see [18]). Then define
$\J_X(a,b)$ to be the union of those components of $\M(X/C)$ which
contain a rank $a$, degree $b$ vector bundle
supported on a non-singular fibre of $\pi$.

The
details of this construction, and the proof that $\J_X(a,b)$ is an
elliptic surface, are given in section 4. The essential point is that each component of
the moduli of stable sheaves on an elliptic
curve is again an elliptic curve, so that $\J_X(a,b)$ has a natural
elliptic fibration structure.

\begin{thm}
\label{moduli}
The moduli space $\M=\M_X(r,\De,k)$ is a smooth (non-empty)
projective variety and is birationally equivalent to
$$\Pic^{\oo}(\J_X(a,b))\times\Hilb^t(\J_X(a,b)),$$
where $(a,b)$
is the unique pair of integers satisfying $br-a(\De\cdot f)=1$ and $0<a<r$.
Furthermore, if $r>at$ the birational equivalence extends to give an
isomorphism of varieties.
\end{thm}

This generalises a result of Friedman who looked at the rank 2 moduli. In that case one always
has $a=1$, and the surface $\J_X(a,b)$ is the relative
Picard scheme $J^b(X)$ of [7].

\subsection{}
To prove Theorem \ref{moduli} we use a relative version of the Fourier-Mukai
transform. Let $\D(S)$ denote the bounded derived category of
coherent sheaves on a Noetherian scheme $S$. In general, given two smooth varieties $X$ and $Y$ and an
object $\PP$ of $\D(X\times Y)$, one defines a functor $\Phi^\PP_{Y\to
X}:\D(Y)\lra \D(X)$, by the formula
$$\Phi^\PP_{Y\to X}(-)=\R\pi_{X,*}(\PP\Ltensor\pi^*_Y(-)),$$
where $\pi_X$ and $\pi_Y$ are the projections from $X\times Y$ to $X$
and $Y$ respectively. In
some very special circumstances this functor is an
equivalence of categories, and one then has a useful tool for
studying sheaves on $X$ and $Y$. Examples include
the original Fourier-Mukai transform [11] (where $X$ is an Abelian variety
and $Y$ is its dual) and the generalised transforms of [4] and [17] (where $X$ and
$Y$ are isogenous K3 surfaces).

\begin{thm}
\label{biggy}
Given integers $a>0$ and $b$ with $a\m_X$ coprime to $b$, let
$Y=\J_X(a,b)$. Then there exist
tautological sheaves on $X\times Y$, supported on $X\times_C Y$, and
for each such sheaf $\PP$,
the functor $\Phi^{\PP}_{Y\to X}:\D(Y)\lra \D(X)$ is an
equivalence of categories.
\end{thm}

Although this result involves derived categories in an essential way,
we shall see that the functors $\Phi^{\PP}_{Y\to X}$ often take sheaves to sheaves, and
thus yield concrete results about moduli of vector bundles. Indeed, once the basic
properties of these functors are known, the proof of Theorem
\ref{moduli} is
relatively simple, and it seems likely that the same functors will prove
useful for solving other moduli problems on elliptic surfaces.
One might also expect similar results for
higher-dimensional elliptic fibrations.

\begin{nota}
All schemes will be Noetherian $\C$-schemes and all morphisms will be morphisms
over $\C$. By a sheaf on a scheme $X$ we mean a coherent
$\OO_X$-module. A point of $X$ will mean a closed point. $\D(X)$
denotes the bounded derived category of sheaves on $X$. We refer to [9] for properties of $\D(X)$.

Given an object $E$ of $\D(X)$ let $\HH^i(E)$ denote the $i$th
cohomology sheaf of $E$. We say that $E$ is a sheaf if
$\HH^i(E)=0$ when $i\neq 0$. The derived
dual, $\R\lHom(E,\OO_X)$, is denoted $E^{\dual}$, and $E[n]$
denotes the object $E$ shifted to the left by $n$ places. 

By a variety we mean an integral, separated scheme of finite type over
$\C$. The canonical bundle of a smooth variety $X$ is
written $\omega_X$, and for objects $E$ and $F$ of $\D(X)$ we define
$$\chi(E,F)=\sum(-1)^i\dim_{\,\C}\Hom_{\D(X)}(E,F[i]),$$
and write $\ch(E)$ for the Chern
character of $E$.
This defines a map
$$\ch:\D(X)\lra H^{2*}(X,\C).$$
We denote the image (an Abelian group) by $\ch(X)$.

We refer to [18] for the definitions of pure-dimension sheaves and
stable sheaves on a projective scheme $X$. They depend on the
polarization chosen for $X$.
\end{nota}

\begin{ack}
This work forms part of my PhD project, which is funded by the EPSRC. I would like to thank my
supervisor, Antony Maciocia, for teaching me about Fourier-Mukai
transforms, and for providing a great deal of help and encouragement.
The original idea of looking for relative transforms on elliptic surfaces was his.
\end{ack}


\section{General Properties of $\Phi^{\PP}_{Y\to X}$}

Throughout this section $X$ and $Y$ are smooth varieties and $\PP$ is an
object of $\D(X\times Y)$. We shall be mainly interested in the case when $\PP$ is a
sheaf on $X\times Y$, flat over both factors. Let $\Phi$ denote the corresponding
functor $\Phi_{Y\to X}^{\PP}:\D(Y)\to\D(X)$ defined in the introduction.
We state various properties of $\Phi$ which we shall need. These all appear in some form in Mukai's original
papers on Abelian varieties and K3 surfaces [11], [13], [14].
The general results have been worked out by A. I. Bondal and D. O. Orlov
[6], [17], and A. Maciocia [10].

\subsection{}
First we define WIT (weak index theorem) sheaves. Given an object $E$ of $\D(Y)$ put
$$\Phi^i(E)=\HH^i(\Phi(E)).$$
A sheaf $E$ on $Y$ is said to be
$\Phi$-WIT$_i$ if $\Phi^j(E)=0$ for all $j\neq i$, or equivalently
if $\Phi(E)[i]$ is a sheaf.
We say $E$ is $\Phi$-WIT if it is
$\Phi$-WIT$_i$ for some $i$, and in this case we often write $\Hat E$ for $\Phi^i(E)$, and refer to
$\Hat E$ as the {\it transform} of $E$.

\subsection{}
By Grothendieck's
Riemann-Roch theorem there is a group homomorphism
$\ch(\Phi)$ making the following diagram commute
$$
\begin{array}{ccc}
\D(Y) &\lRa{\Phi} &\D(X) \\
\scriptstyle{\ch}\bda  &&\scriptstyle{\ch}\bda \\
\ch(Y) &\lRa{\ch(\Phi)} &\ch(X)
\end{array}
$$
It is given by
$$\ch(\Phi)(y)=\pi_{X,*}(p\cdot\pi_Y^*y),$$
where $p=\ch(\PP)\cdot\pi_Y^*(\td_Y)$ and $\td_Y$ is the Todd class of
$Y$.

\subsection{}
Define
$$\QQ=(\PP^{\dual}\tensor\pi_X^*\omega_X)[\dim X+\dim Y-\dim \PP],$$
and put $\Psi=\Phi^{\QQ}_{X\to Y}$.
It is a simple consequence of Grothendieck-Verdier
duality (see e.g. [6], Lemma 1.2) that $\Psi[\dim \PP-\dim Y]$ is a
left adjoint of $\Phi$.

If $\Phi$ is fully faithful one has an isomorphism of functors
$$\Psi\circ\Phi\isom \id_{\D(Y)}[\dim Y-\dim\PP].$$

The reason for the apparently strange choice of shift in the definition of $\QQ$ is to give
it a good chance of being a sheaf. For example when $\PP$ is a vector
bundle, so is $\QQ$.

\subsection{}
Assume that $\PP$ is a sheaf on $X\times Y$, flat over both factors. Then
$\Phi$ preserves families of sheaves. In detail, let $S$ be a scheme
and $\EE$ an $S$-flat sheaf on $Y\times S$. Then
$$U=\{s\in S:\EE_s\mbox{ is } \Phi\mbox{-WIT}_i\}$$
is the set of points of an open subscheme of $S$. Furthermore there
exists a $U$-flat sheaf $\FF$ on $X\times U$, such that for all $s\in U$, $\FF_s=\Phi^i(\EE_s)$. The proof
of this result is identical to that of [14], Theorem 1.6.

\subsection{}
Suppose now that $\Phi$ is fully faithful, and take sheaves $A$ and $B$ on
$Y$, with $A$ $\Phi$-WIT$_a$ and $B$ $\Phi$-WIT$_b$. Then for all $i$,
one has
$$\Hom_{\D(Y)}(A,B[i])=\Hom_{\D(X)}(\Hat A[-a],\Hat B[i-b]).$$
Rewriting this gives the identity
\begin{equation}
\label{percy}
\Ext_Y^i(A,B)=\Ext_X^{i+a-b}(\Hat A,\Hat B),
\end{equation}
which is referred to as the Parseval theorem. As a special case, note
that if $A$ is simple then so is $\Hat A$.

\subsection{}
Assume that $X$ and $Y$ have the same dimension. Recall that a $Y$-flat
sheaf $\PP$ on $X\times Y$ is said to be {\it strongly simple}
over $Y$ if $\PP_y$ is simple for all $y\in Y$, and if for any pair
$y_1,y_2$ of distinct points of $Y$ and any integer $i$ one has
$\Ext^i_X(\PP_{y_1},\PP_{y_2})=0$.

\begin{thm}
\label{equivalence}
Let $\PP$ be a $Y$-flat sheaf on $X\times Y$. Then
$\Phi$ is fully faithful iff $\PP$ is strongly simple over
$Y$. If $\PP$ is flat over $X$ and $Y$ then $\Phi$ is an equivalence iff
$\PP$ is strongly simple over both factors.
\qed
\end{thm}

The main idea behind the proof was given by Mukai ([13],
Theorem 4.9) and the
result has appeared in
various forms since then. The most general statement is due to Bondal
and Orlov ([6], Theorem 1.1).

Following Orlov we shall make essential use of the following lemma. It
is an immediate consequence of [5], Lemma 3.1.

\begin{lemma}
\label{useful}
Suppose $\Phi$ is fully faithful. Then $\Phi$ is an equivalence iff
for any object $E$ of $\D(X)$,
$\Psi(E)\isom 0$ implies $E\isom 0$.
\qed
\end{lemma}


\section{Fourier-Mukai Transforms for Elliptic Curves}

Here we illustrate the results of the last section by considering the
case when $X$ and $Y$ are elliptic curves. 

\subsection{}
Let $X$ be an elliptic curve. Given a sheaf $E$ on $X$ we write
its Chern class as a pair of integers $(r(E),d(E))$. Let $a$ and $b$ be coprime integers with
$a>0$ and let $Y$ be the moduli space of stable bundles on $X$ of
Chern class $(a,b)$. In fact it is a consequence of the work of Atiyah ([2], Theorem 7) that $Y$ is isomorphic to $X$; we
preserve the distinction for clarity.
Let $\PP$ be a tautological bundle on $X\times Y$, and put
$$\Phi=\Phi^{\PP}_{Y\to X},\qquad \Psi=\Phi^{\PP^{\dual}}_{X\to Y}.$$
As we noted in 2.3, $\Psi[1]$ is a left adjoint of $\Phi$.

\begin{prop}
\label{yawn}
The functor $\Phi$ is an
equivalence.
\end{prop}

\begin{pf}
First note that $\PP$ is strongly simple over $Y$, since for any pair $P_1,P_2$ of
non-isomorphic stable bundles on $X$ with the same Chern class,
Serre duality gives
$$\Ext^1_X(P_2,P_1)=\Hom_X(P_1,P_2)^{\dual}=0.$$
It follows from Theorem
\ref{equivalence} that $\Phi$ is fully faithful, so
$\Psi\circ\Phi\isom \id_{\D(Y)}[-1]$.

Next
observe that the
group homomorphism $\ch(\Phi)$ must be an isomorphism, since
$$\ch(X)\isom\ch(Y)=\Z\oplus\Z,$$
and $\ch(\Psi)\circ\ch(\Phi)=-\id_{\ch(Y)}$.

To complete the
proof use Lemma \ref{useful}. Suppose $E$ is an object
of $\D(X)$ such that $\Psi(E)\isom 0$. Consider the hypercohomology spectral sequence
$$E^{p,q}_2=\Psi^p(\HH^q(E))\implies \Psi^{p+q}(E)=0.$$
Since $E^{p,q}_2=0$ unless $p=0$ or 1, the spectral sequence
degenerates at the $r=2$ level. It follows that $\Psi(\HH^q(E))=0$ for
all $q$. But then, since $\ch(\Psi)$ is an isomorphism, one has that
$\HH^q(E)=0$ for all $q$, so $E\isom 0$.
\qed
\end{pf}

The group homomorphism $\ch(\Phi)$ is invertible and takes $(0,1)$ to $(a,b)$, so must be given by some matrix
$$\mat{c}{a}{d}{b}$$
where $c$ and $d$ are integers satisfying $bc-ad=\pm 1$. Then
$\ch(\Psi)$, which is the inverse of $-\ch(\Phi)$, takes $(0,1)$ to
$\pm (a,-c)$. Since $\Psi$ is given by a sheaf on $X\times Y$, we
must take the positive sign, so that in fact $bc-ad=1$.

This condition
does not define $c$ and
$d$ uniquely: we may replace them by $c+na$ and
$d+nb$ for any integer $n$. This corresponds to twisting $\PP$ by the pull-back of a line bundle of degree
$n$ on $Y$. By varying $n$ we obtain all possible values of $c$ and
$d$.

\begin{thm}
\label{ellcur}
Let $X$ be an elliptic curve and take an element
$$A=\mat{c}{a}{d}{b}\in\mbox{SL}_2(\Z),$$
such that $a>0$. Then there exist vector bundles
on $X\times X$ which are strongly simple over both factors, and
which restrict to give bundles of Chern class $(a,c)$ on the first factor and
$(a,b)$ on the second. For any such bundle $\PP$, the resulting functor $\Phi=\Phi^{\PP}_{X\to
X}$ is an equivalence, and satisfies
$$\col{r(\Phi E)}{d(\Phi E)}=\mat{c}{a}{d}{b}\col{r(E)}{d(E)},$$
for all objects $E$ of $\D(X)$.
\qed
\end{thm}

For the usual Fourier-Mukai transform $\F$ on $X$ (see [11]) one has
$$A=\mat{0}{1}{-1}{0}.$$
The equivalences we have found are not
essentially new, since one can check that they
can all be
obtained from composites of the functors $\F$ and $L\tensor(-)$ for
line bundles $L$ on $X$. Later, however, we shall
try to apply the transforms on each fibre of an elliptic
surface, and this will only be possible for certain choices of $A$. 

\subsection{}
We conclude by showing that the transforms take simple
sheaves to simple sheaves.

\begin{prop}
\label{stable}
Let $X$ be an elliptic curve and $\PP$ a bundle on
$X\times X$, strongly simple over one factor. Then
$\Phi=\Phi^{\PP}_{X\to X}$ is an equivalence. Furthermore any simple
sheaf $E$ on $X$ is $\Phi$-WIT and the transform $\Hat E$ is a simple sheaf.
\end{prop}

\begin{pf}
The argument of Proposition \ref{yawn} shows that $\Phi$ is an equivalence, so defining
$\Psi$ as above, there is an isomorphism
$$\Psi\circ\Phi\isom \id_{\D(X)}[-1],$$
and hence, for any sheaf $E$ on $X$, a spectral sequence
$$E^{p,q}_2=\Psi^p(\Phi^q(E))\implies\left\{\begin{array}{ll}
E &\mbox{ if $p+q=1$} \\
0 &\mbox{ otherwise.}
\end{array} \right. $$
Now $E^{p,q}_2=0$ unless $0\leq p,q\leq 1$, so this gives a short
exact sequence
$$0\lra \Psi^1(\Phi^0(E))\lra E\lra \Psi^0(\Phi^1(E))\lra 0,$$
together with the information that $\Phi^0(E)$ is $\Psi$-WIT$_1$ and
$\Phi^1(E)$ is $\Psi$-WIT$_0$. The Parseval theorem then implies that
$$\Ext^1_X(\Psi^0(\Phi^1(E)),\Psi^1(\Phi^0(E)))=0,$$
so $E$ is given by a trivial extension. If $E$ is simple it follows that one
of the two sheaves $\Psi^1(\Phi^0(E))$ or $\Psi^0(\Phi^1(E))$ is zero, and $E$ is $\Phi$-WIT. The transform
$\Hat E$ is then simple, as we noted in 2.5.
\qed
\end{pf}

\begin{remark}
A straightforward
application of Serre duality shows that a simple sheaf on $X$ is
either a stable vector bundle or the skyscraper sheaf of a point of
$X$. In particular, a vector bundle on $X$ is simple iff it is stable.
\end{remark}


\section{The Elliptic Surfaces $\J_X(a,b)$}

In this section we introduce the elliptic surfaces $\J_X(a,b)$
mentioned in the introduction.  By an elliptic surface we shall mean
a smooth variety $X$ of dimension 2 together with a smooth curve $C$
and a relatively minimal morphism $\pi:X\to C$ whose general fibre is
an elliptic curve. We often abuse notation and refer to $X$ as an
elliptic surface, or an elliptic surface over $C$, and take the
morphism $\pi$ as given.

\subsection{}
Let $\ellsur$ be an elliptic surface. Recall ([3], V.12.3) that the
canonical bundle of $X$ takes the form
$$\omega_X=\pi^*(L)\tensor \OO_X(\sum{(m_i-1)f_i)},$$
where $L$ is a line bundle on $C$ and $m_1f_1,\cdots,m_kf_k$ are the
multiple fibres of $\pi$. This formula depends on the assumption that
$\pi$ is relatively minimal.

We denote the algebraic equivalence class of a
fibre of $\pi$ by $f$, and for any object $E$ of $\D(X)$ define the {\it fibre
degree} of $E$ to be
$$d(E)=c_1(E)\cdot f.$$
Note that the restriction of a sheaf $E$ on $X$ to a general fibre of $\pi$ has Chern
class $(r(E),d(E))$. We say that a sheaf $E$ on $X$ is a {\it fibre sheaf} if
$r(E)=d(E)=0$, or equivalently if the support of $E$ is contained in
the union of finitely many fibres of $\pi$. In this case $E\tensor\omega_X$ has the same Chern
class as $E$, and for any other sheaf $F$ on $X$,
\begin{equation}
\label{serre}
\chi(E,F)=\chi(F,E).
\end{equation}

Let $\m_X$ denote the highest common factor of the fibre degrees of sheaves
on $X$. Equivalently $\m_X$ is the smallest positive integer such that
there is a divisor $\sigma$ on $X$ with $\sigma\cdot
f=\m_X$. Note that, by Riemann-Roch, given a divisor of positive fibre
degree, we can add a large multiple of $f$ and obtain an effective
divisor of the same fibre degree.

\subsection{}
Let $\ellsur$ be an elliptic surface and fix integers $a>0$ and $b$,
with $a\m_X$ coprime to $b$. Take a polarization of $X$ of fibre
degree coprime to $b$. By the results of [18], there exists a
relative moduli scheme $\M(X/C)\to C$, whose points represent stable
pure-dimension 1 sheaves on fibres of $\pi$.

\begin{dfn}
Let $\J_X(a,b)$ be the union of those components
of $\M(X/C)$ which contain a point representing a rank $a$, degree $b$
vector bundle on a
non-singular fibre of $\pi$. Let $\Hat{\pi}$ denote the natural map
$\Hat{\pi}:\J_X(a,b)\lra C.$
\end{dfn}

First note that the coprimality assumptions we have made imply that
$Y=\J_X(a,b)$ is a fine moduli scheme. Thus $Y$
is a projective scheme whose points all represent strictly stable sheaves. Furthermore, by an argument of Mukai
([13], Theorem A.6),
there is a tautological sheaf $\PP$ on $X\times_C Y$, such that for each point $y\in Y$, the stable sheaf corresponding to $y$
is given by $\PP_y$, the restriction of $\PP$ to
$X_{\Hat{\pi}(y)}\times\{y\}$.

Let $U$ be the set of points $p\in C$ such that the fibre $X_p$
is non-singular. The fibre of
$\Hat{\pi}$ over a point $p\in U$ is the moduli space of rank $a$,
degree $b$ stable sheaves
on $X_p$, which, as we noted in section 3, is isomorphic to $X_p$. Thus $\Hat{\pi}$ is an elliptic fibration. Clearly $\Hat{\pi}$
is dominant, hence surjective, so there is some component of $Y$ which
contains sheaves supported on every fibre of $\pi$. Any other component of
$Y$ must contain a sheaf supported on a non-singular fibre, but the fibre of
$\Hat{\pi}$ over every point of $U$ is connected. It follows
that $Y$ is
connected.

\smallskip

Now let $\M(X)$ denote the moduli space of stable pure-dimension 1
sheaves on $X$, and let $Z$ be the union of those components of $\M(X)$
which contain a rank $a$, degree $b$ sheaf supported on a non-singular fibre of
$\pi$. Thus points of $Z$ correspond to strictly stable sheaves of
Chern class $(0,af,-b)$.

There is a natural `extension by zero' morphism
$i:Y\lra Z$, which maps a point $y\in Y$ representing the stable
sheaf $\PP_y$ on the fibre $X_{\Hat{\pi}(y)}$, to the point $z\in Z$
representing the stable sheaf on
$X$ obtained by extending $\PP_y$ by zero. This morphism $i$ induces
a bijection on points, since every stable sheaf of Chern class
$(0,af,-b)$ is supported on some fibre of $\pi$ (the support of a
stable sheaf must be connected). I claim
that $Z$ is a non-singular projective surface; it will follow from
this that $i$ is an isomorphism and that $Y$ is an elliptic surface
over $C$.

\smallskip

Given a point $y\in Y$ we shall identify the sheaf $\PP_y$ with its
extension by zero on $X$. If $y\in Y$ is such that $\PP_y$ is supported on a non-singular fibre
of $\pi$, then
\begin{equation}
\label{phew}
\PP_y=\PP_y\tensor \omega_X,
\end{equation}
because the restriction of $\omega_X$ to any non-singular fibre of
$\pi$ is trivial. By EGA III.7.7.8 the
dimension of the space
$$\Hom_X(\PP_y,\PP_y\tensor \omega_X)$$
is upper semi-continuous on $Y$, so for all $y\in Y$ there is a non-zero
morphism $\PP_y\to\PP_y\tensor\omega_X$. But both these sheaves are
stable with the same Chern class, so they are isomorphic and \ref{phew} holds for all $y\in Y$.

The Riemann-Roch formula gives
$$\chi(\PP_y,\PP_y)=-(af)^2=0,$$
so the Zariski tangent space to $Z$ at a point $i(y)$, which is given by
$$\Ext^1_X(\PP_y,\PP_y)$$
(see e.g. [19]), always has dimension 2. Now $Y$ fibres over $C$
with elliptic fibres, so has dimension at least 2, and it then follows that $Z$ is a non-singular
projective surface as claimed.

\smallskip

Extending our tautological sheaf $\PP$ by zero, we obtain a sheaf on
$X\times Y$ which we shall also denote by $\PP$, such that for each point $y\in Y$, $\PP_y$ is a stable
sheaf of Chern class $(0,af,-b)$ on $X$.

For any two distinct points $y_1$, $y_2$ of $Y$, Serre duality implies that
$$\Ext_X^2(\PP_{y_1},\PP_{y_2})=\Hom_X(\PP_{y_2},\PP_{y_1})^{\dual}=0,$$
and since $\chi(\PP_{y_1},\PP_{y_2})=0$,
this is enough to show that $\PP$ is strongly simple over $Y$. By
Theorem \ref{equivalence}, the functor $\Phi=\Phi^{\PP}_{Y\to X}$ is
fully faithful.

\begin{prop}
\label{nearly}
The scheme $Y=\J_X(a,b)$ is an elliptic surface over $C$. Furthermore,
the sheaf $\PP$ is strongly simple over $Y$, so $\Phi=\Phi^{\PP}_{Y\to
X}$ is fully faithful.
\end{prop}

\begin{pf}
It only remains to show that $Y$ is relatively minimal over $C$. Suppose
not, i.e. that there exists a $(-1)$-curve $D$ contained in a fibre
of $\Hat{\pi}$. Then
$\KK_Y\cdot D<0$, so
$$\chi(\OO_D,\OO_Y)=\chi(\omega_Y|_D)\neq\chi(\OO_D)=\chi(\OO_Y,\OO_D).$$
Since $\Phi$ is fully faithful this implies that $\chi(E,F)\neq\chi(F,E),$
where $E=\Phi(\OO_D)$ and $F=\Phi(\OO_Y)$. But for each $i$,
$\HH^i(E)$ is a fibre sheaf (because $\OO_D$ is), so this contradicts
\ref{serre}.
\qed
\end{pf}

\begin{remark}
\label{tiup}
Suppose we use two different polarizations of $X$ to define
elliptic surfaces  $\J_X(a,b)$ and $\J'_X(a,b)$ over $C$. Then, since the stability of a sheaf on a smooth curve does not
depend on a choice of polarization, the
two spaces will be isomorphic 
over the open subset $U$ considered above, and hence birational. Since
both are relatively minimal over $C$, [3], Proposition III.8.4 implies
that they are isomorphic as elliptic surfaces over $C$.
\end{remark}

\subsection{}
In the next section we show that the functor $\Phi$ is an
equivalence. For now, let us note that the restriction of $\Phi$ to a
non-singular fibre of $\pi$ yields one of the transforms considered in
section 3. Indeed, if $p\in C$ and $i_p:X_p\embed X$ and $j_p:Y_p\embed Y$ are the
inclusion of the non-singular fibres $X_p$ and $Y_p$, then a simple
base-change (see [6], Lemma 1.3) gives an isomorphism of functors
$$\L i_p^*\circ\Phi\isom\Phi_p\circ\L j_p^*.$$
Here $\Phi_p$ is the functor $\Phi_{Y_p\to X_p}^{\PP_p}$ and
$\PP_p$, the restriction of $\PP$ to $X_p\times Y_p$, is a
tautological bundle parameterising stable bundles on $X_p$ of rank $a$
and degree $b$. Thus $\Phi_p$ coincides with one of the transforms of
Theorem \ref{ellcur}, and in particular there is a matrix
$$\mat{c}{a}{d}{b}\in\mbox{SL}_2(\Z),$$
such that for all objects $E$ of $\D(Y)$,
\begin{equation}
\label{blag}
\col{r(\Phi E)}{d(\Phi E)}=\mat{c}{a}{d}{b}\col{r(E)}{d(E)}.
\end{equation}
Furthermore, Proposition \ref{stable} gives

\begin{lemma}
\label{whoop}
Let $E$ be a $\Phi$-WIT sheaf on $Y$ whose restriction
to the general fibre of $\Hat{\pi}$ is simple. Then the restriction of $\Hat
E$ to the general fibre of $\pi$ is also simple.
\qed
\end{lemma}


\section{Fourier-Mukai Transforms for Elliptic Surfaces}

Here we prove Theorem \ref{biggy}. Let $\ellsur$ be an elliptic surface, fix
integers $a>0$ and $b$, with $a\m_X$ coprime to $b$, and let $Y$ denote
the elliptic surface
$\Hat{\pi}:\J_X(a,b)\lra C$ defined in the last section.
Fix a tautological sheaf on $X\times_C Y$, and extend by zero to
obtain a sheaf $\PP$ on $X\times Y$. Let $\QQ$ be the object
$(\PP^{\dual}\tensor\pi_X^*\omega_X)[1]$ of $\D(X\times Y)$, and define functors
$$\Phi=\Phi^{\PP}_{Y\to X},\qquad\Psi=\Phi^{\QQ}_{X\to Y}.$$
As we noted in 2.3, $\Psi[1]$ is a left adjoint of $\Phi$.

\begin{lemma}
\label{tech}
The object $\QQ$
is a sheaf on $X\times Y$. Moreover, $\PP$ and $\QQ$ are both flat
over $X$ and $Y$.
\end{lemma}

\begin{pf}
For each point $(x,y)\in X\times Y$, consider the commutative diagram
$$
\begin{array}{ccc}
\Spec\C & \lRa{j_x} & X \\
\scriptstyle{j_y}\bda  && \scriptstyle{i_y}\bda \\
Y & \lRa{i_x} & X\times Y
\end{array}
$$
where $j_x$ and $i_x$ are the inclusions of $\{x\}$ in $X$ and
$\{x\}\times Y$ in $X\times Y$ respectively. Similarly for $j_y$ and $i_y$.

First note that by an argument of Mukai ([12], p. 105), any
pure-dimension 1 sheaf on any surface has a locally-free
resolution of length 2. This implies that for all $y\in Y$, $(\PP_y)^{\dual}[1]$ is a sheaf on $X$. But
$$\L i_y^*(\PP^{\dual}[1])=(\L i_y^*(\PP))^{\dual}[1]=(\PP_y)^{\dual}[1],$$
so the corresponding hypercohomology spectral sequence implies that
$\PP^{\dual}[1]$ is a sheaf on $X\times Y$, flat over $Y$ (see [6],
Proposition 1.5).

To show that $\PP$ is flat over $X$ consider the spectral sequence
$$E_{p,q}^2=\L_p j_y^*(\L_q
i_x^*(\PP))\implies\L_{(p+q)}j_x^*(\PP_y).$$
Here $\L_p f^*(E)$ denotes the $(-p)$th cohomology sheaf of $\L f^*(E)$.
Since $\PP_y$ has a two-step resolution, the right-hand side is non-zero only if
$p+q=0$ or 1, so one
concludes that
$$\L_1 j_y^*(\L_1 i_x^*(\PP))=0,$$
for all $y\in Y$. This implies that $\L_1 i_x^*(\PP)$ is
locally free on $Y$. But for any $x\in X$
one can find $y\in Y$ such that $(x,y)$ does not lie in the support of
$\PP$, so $\L_1 i_x^*(\PP)=0$ for all $x\in X$, and $\PP$ is flat over $X$.

Finally, the isomorphism
$$\L i_x^*(\PP^{\dual}[1])\isom(\PP_x)^{\dual}[1],$$
implies that both sides are sheaves on $Y$, so
$\PP^{\dual}[1]$ is flat over $X$.
\qed
\end{pf}

The next lemma shows that the relationship between $X$
and $Y$ is
entirely symmetrical.

\begin{lemma}
\label{crafty}
There exists an integer $c$ such that $X\isom \J_Y(a,c).$
\end{lemma}

\begin{pf}
If $X_p$ is a non-singular fibre of $\pi$ then
the restriction of $\PP$ to $X_p\times Y_p$ is a tautological bundle
parameterising stable bundles on $X_p$. By the results of section 3,
this bundle is strongly simple over both factors, so for any point $x\in X$
lying on a non-singular fibre of $\pi$, the sheaf $\PP_x$ is a stable
sheaf on $Y$. Let its Chern class be $(0,af,-c)$. I claim that $c$ is coprime to
$a\m_Y$, so that $\J_Y(a,c)$ is well-defined. Assuming this for the
moment, note that as in Remark \ref{tiup}, the two elliptic
surfaces $X$ and $\J_Y(a,c)$ over $C$ are isomorphic away from the
singular fibres, so are isomorphic.

Since the object $\QQ_x$ of $\D(Y)$ has Chern class $(0,af,c)$, to prove the claim it will be enough to exhibit an object $E$ of
$\D(Y)$ such that $\chi(\QQ_x,E)=1$. But this is possible by
\ref{blag}, since the result of section 2.3 implies that
$$\chi(\QQ_x,E)=-\chi(\C_x,\Phi E)=-r(\Phi E),$$
for any object $E$ of $\D(Y)$.
\qed
\end{pf}

\smallskip

We can now prove Theorem \ref{biggy}.
By Proposition \ref{nearly}, $\Phi$ is fully faithful, so
$\Psi\circ\Phi\isom\id_{\D(Y)}[-1]$. It follows that $\ch(\Psi)\circ\ch(\Phi)=-\id_{\ch(Y)},$
and $\ch(\Phi)$ embeds
$\ch(Y)$ as a direct summand of $\ch(X)$. Applying Lemma \ref{crafty},
we can repeat the argument and obtain $\ch(X)$ as a direct summand of
$\ch(Y)$. This shows that $\ch(\Phi)$ is an isomorphism.

We complete the proof by applying Lemma \ref{useful}. Suppose $E$ is
an object of $\D(X)$ such that $\Psi(E)\isom 0$.
Consider the hypercohomology spectral sequence
$$E^{p,q}_2=\Psi^p(\HH^q(E))\implies\Psi^{p+q}(E)=0.$$
Since $\QQ$ is supported on
$X\times_C Y$ and is flat over
$X$ and $Y$, one has that $E^{p,q}_2=0$ unless $p=0$ or 1. It
follows that the spectral sequence degenerates at the $r=2$ level, so
$\Psi(\HH^q(E))=0$ for all $q$. But since $\ch(\Psi)$ is an
isomorphism, and a non-zero sheaf has non-zero Chern character, this
implies that
$\HH^q(E)=0$ for all $q$, so $E\isom 0$.

\smallskip

We summarise our results in the following theorem.

\begin{thm}
\label{super}
Let $\ellsur$ be an elliptic surface and
take an element 
$$\mat{c}{a}{d}{b}\in\mbox{SL}_2(\Z),$$
such that $\m_X$ divides $d$ and $a>0$. Let $Y$ be the elliptic
surface $\J_X(a,b)$ over $C$. Then there exist sheaves $\PP$ on
$X\times Y$, flat and strongly simple over both factors such that for
any point $(x,y)\in X\times Y$,
$\PP_y$ has Chern class $(0,af,-b)$ on $X$ and $\PP_x$ has Chern class
$(0,af,-c)$ on $Y$.

For any such sheaf $\PP$, the resulting functor $\Phi=\Phi^{\PP}_{Y\to X}$ is an
equivalence and satisfies
\begin{equation}
\label{num}
\col{r(\Phi E)}{d(\Phi E)}=\mat{c}{a}{d}{b}\col{r(E)}{d(E)},
\end{equation}
for all objects $E$ of $\D(Y)$.
\end{thm}

\begin{pf}
Take a tautological sheaf $\PP$ on $X\times Y$ and put
$\Phi=\Phi^{\PP}_{Y\to X}$. As we showed above, $\Phi$ is an
equivalence and there
exist integers $c$ and $d$ such that \ref{num} holds. Now $\m_X$
divides $d(\Phi E)$ for any object
$E$ of $\D(Y)$, so $\m_X$ divides $\m_Y$ and $d$. By symmetry
$\m_X=\m_Y$. As in section 3, $c$ and $d$ are not uniquely defined: we can
replace them by $c+n\m_Xa$ and
$d+n\m_Xb$ by twisting $\PP$ by the pull-back of a line bundle of
fibre degree $n\m_X$ on $Y$.
\qed
\end{pf}

\begin{remark}
\label{stability}
As a corollary of the proof of Theorem \ref{super}, note that we can
always choose $\PP$ so that for some
polarization of $X$, $\PP_y$ is stable for all $y\in Y$.
By Lemma \ref{crafty}, we could also view $X$ as a moduli space of
sheaves on $Y$, and take $\PP$ such that for some polarization of $Y$, $\PP_x^{\dual}[1]$ is stable for all $x\in X$.
\end{remark}


\section{Properties of the Transforms}

Let $\ellsur$ be an elliptic surface, fix an element
$$\mat{c}{a}{d}{b}\in\mbox{SL}_2(\Z),$$
with $a>0$ and $\m_X$ dividing $d$, let $Y$ be the elliptic surface $\J_X(a,b)$ and take a sheaf
$\PP$ on $X\times Y$ as in Theorem \ref{super}. As in the last
section, define the sheaf $\QQ=(\PP^{\dual}\tensor\pi_X^*\omega_X)[1],$
and the functors $\Phi$ and $\Psi$.
Since $\Phi$ is an equivalence, one has isomorphisms
\begin{equation}
\label{star}
\Psi\circ\Phi\isom \id_{\D(Y)}[-1],\qquad\Phi\circ\Psi\isom
\id_{\D(X)}[-1].
\end{equation}
In this section we give some properties of the transforms which will
be useful in section 7. 
Note that, because of the symmetry of the situation, for each result
we give here, there will be another result obtained by exchanging
$\Phi$ and $\Psi$, and $X$ and $Y$.

\subsection{}
The functor $\Phi$ is left exact,
because $\PP$ is flat over $Y$. Thus
given a short exact sequence
$$0\lra A\lra B\lra C\lra 0,$$
one obtains a long exact sequence
\begin{eqnarray}
\label{long}
 0&\lra&\Phi^0(A)\lra\Phi^0(B)\lra\Phi^0(C) \\
  &\lra&\Phi^1(A)\lra\Phi^1(B)\lra\Phi^1(C)\lra0. \nonumber
\end{eqnarray}

\subsection{}
The isomorphisms \ref{star} imply that if $E$ is a $\Psi$-WIT$_i$
sheaf on $X$ ($i=0,1$), then $\Hat E$ is a $\Phi$-WIT$_{1-i}$ sheaf on $Y$. More
generally, for any sheaf $E$
on $X$ there is a spectral sequence
$$E^{p,q}_2=\Phi^p(\Psi^q(E))\implies\left\{\begin{array}{ll}
E &\mbox{ if $p+q=1$} \\
0 &\mbox{ otherwise.}
\end{array} \right. $$
Since $E^{p,q}_2=0$ unless $0\leq p,q\leq 1$, this
yields a short exact sequence
$$0\lra \Phi^1(\Psi^0(E))\lra E\lra \Phi^0(\Psi^1(E))\lra 0,$$
together with the information
that $\Psi^0(E)$ is $\Phi$-WIT$_1$ and $\Psi^1(E)$ is $\Phi$-WIT$_0$.
Note that \ref{star} then implies that $\Phi^1(\Psi^0(E))$ is
$\Psi$-WIT$_0$ and $\Phi^0(\Psi^1(E))$ is $\Psi$-WIT$_1$.

\begin{lemma}
\label{specseq}
For any sheaf $E$ on $X$ there is a unique short exact
sequence
$$0\lra A\lra E\lra B\lra 0,$$
such that $A$ is $\Psi$-WIT$_0$ and $B$ is $\Psi$-WIT$_1$.
\end{lemma}

\begin{pf}
For the uniqueness, suppose there is another such sequence
$$0\lra A'\lra E\lra B'\lra 0.$$
Then, since $A'$ is $\Psi$-WIT$_0$ and $B$ is $\Psi$-WIT$_1$, the
Parseval theorem implies that there is
no non-zero map $A'\to B$, so the inclusion of $A'$ in $E$ factors
through $A$. By symmetry $A=A'$.
\qed
\end{pf}

\subsection{}
Given a torsion-free sheaf $E$ on $X$, put $\mu(E)=d(E)/r(E)$.

\begin{lemma}
\label{chern}
Let $E$ be a torsion-free sheaf on $X$. If $E$ is $\Psi$-WIT$_0$ then
$\mu(E)\geq b/a$. Similarly if $E$ is $\Psi$-WIT$_1$ then
$\mu(E)\leq b/a$.
\end{lemma}

\begin{pf}
If $E$ is $\Psi$-WIT$_1$ then $\Psi(E)[1]$ is a sheaf so
$r(\Psi E)\leq 0$. Similarly, if $E$ is $\Psi$-WIT$_0$ then $r(\Psi E)\geq 0$.
Since $\Psi[1]$ is the inverse of $\Phi$, one has
$$\col{r(\Psi E)}{d(\Psi E)}=\mat{-b}{a}{d}{-c}\col{r(E)}{d(E)},$$
for any object $E$ of $\D(X)$. The result follows.
\qed
\end{pf}
A similar argument gives
\begin{lemma}
\label{fibre}
Let $T$ be a $\Psi$-WIT$_1$ torsion sheaf on $X$. Then $T$ is a fibre
sheaf.
\qed
\end{lemma}

Combining Lemma \ref{chern} with Lemma \ref{specseq} we obtain

\begin{lemma}
\label{wit}
Let $E$ be a torsion-free sheaf on $X$ such that the restriction
of $E$ to the general fibre of $\pi$ is stable. Suppose $\mu(E)<b/a$.
Then $E$ is $\Psi$-WIT$_1$.
\end{lemma}

\begin{pf}
Consider the short exact sequence of Lemma \ref{specseq}. If $A$
is non-zero, it is torsion-free and one has $\mu(A)\geq b/a>\mu(E)$. Restricting to the general fibre of $\pi$ this
gives a contradiction. Hence $A=0$ and $E$ is $\Psi$-WIT$_1$.
\qed
\end{pf}

\subsection{}
The final result we shall need is

\begin{lemma}
\label{new}
A sheaf $F$ on $Y$ is $\Phi$-WIT$_0$ iff
$$\Hom_Y(F,\QQ_x)=0\qquad\forall x\in X.$$
\end{lemma}

\begin{pf}
First note that $\QQ_x=\Psi(\C_x)$ is $\Phi$-WIT$_1$. 
If $F$ is $\Phi$-WIT$_0$, then the
Parseval theorem implies that there are no non-zero maps
$F\to\QQ_x$.

Conversely, if $F$ is not $\Phi$-WIT$_0$, then by the argument of Lemma
\ref{specseq}, there is a surjection $F\to B$ with $B$ a
$\Phi$-WIT$_1$ sheaf.
Applying the Parseval theorem again gives
$$\Hom_Y(B,\QQ_x)=\Hom_X(\Hat B,\C_x).$$
Since $\Hat B$ is non-zero, there exists an $x\in X$ and a
non-zero map $B\to \QQ_x$, hence a non-zero map $F\to \QQ_x$.
\qed
\end{pf}


\section{Application to Moduli of Stable Sheaves}

In this section we use the relative transforms we have developed to prove Theorem
\ref{moduli}.

\subsection{}
Let $\ellsur$ be an elliptic surface and fix a triple
$$(r,\De,k)\in\N\times\NS(X)\times\Z,$$
such that $r$ is coprime to $d=\De\cdot f$.
The proof of the following result is entirely analagous to the rank
$2$ case ([7], Theorem I.3.3) so we omit it (see also [15],
Proposition I.1.6).

\begin{prop}
There exist polarizations of $X$ with respect to which a torsion-free
sheaf $E$ on $X$ with Chern class $(r,\De,k)$ is $\mu$-stable whenever
it is $\mu$-semi-stable, and this is the case iff the
restriction of $E$ to all but finitely many fibres of $\pi$ is stable.
\qed
\end{prop}

Taking such a polarization, define
$\M=\M_X(r,\Lambda,k)$ to be the (fine) moduli space of stable torsion-free sheaves
on $X$ of Chern class $(r,\De,k)$. We identify the closed points
of $\M$ with the stable sheaves which they represent. As in the rank
$2$ case ([7], Lemma III.3.6), one shows
that for any $E\in\M$,
$$\Ext^2_X(E,E)=H^2(X,\OO_X).$$
It then follows from the general results of [1] that $\M$, if non-empty, is smooth of dimension $\dim (\Pic^{\oo}(X))+2t$, where
$$2t=2rk-(r-1)\De^2-(r^2-1)\chi(\OO_X),$$
and that if $t<0$ then $\M$ is empty.
In what follows we take $r>1$ and assume that $t$ is non-negative.

Let $a$ and $b$ be the unique pair of
integers satisfying $br-ad=1$,
with $0<a<r$. Let $\Hat{\pi}:Y\to C$ be the elliptic surface $\J_X(a,b)$ and put
$$\NN=\M_Y(1,0,t)=\Pic^{\oo}(Y)\times\Hilb^t(Y).$$
We shall prove Theorem \ref{moduli} by showing that $\M$ is
birationally equivalent to $\NN$.

\subsection{}
Let $\PP$ be a sheaf on $X\times Y$ as in
Theorem \ref{super}, with matrix
$$\mat{r}{a}{d}{b},$$
and define equivalences of categories $\Phi$ and $\Psi$ as in section
6. As we noted in Remark \ref{stability}, we can assume that we have
chosen $\PP$, and a polarization of $Y$, so that $\QQ_x$
is a stable sheaf for all $x\in X$.

Take a sheaf $E$ on $X$ of Chern class $(r,\De,k)$. The formula given in the proof of Lemma \ref{chern} shows that
$\Psi(E)$ has rank 1 and fibre degree 0. Twisting $\PP$ by
the pull-back of a line
bundle on $Y$ we can assume that $c_1(\Psi E)=0$, and the formula
\ref{percy}, together with Riemann-Roch then implies that $\Psi(E)$
has Chern class $(1,0,t)$.

Note that, by Lemma
\ref{wit}, any element $E$ of $\M$ is $\Psi$-WIT$_1$. Define
$$\U=\{E\in\M:\Hat E\mbox{ is torsion-free}\}.$$
By the result of section 2.4, $\U$ is an open subscheme of $\M$. Also define the open subscheme
$$\V=\{F\in \NN:F\mbox{ is }\Phi\mbox{-WIT}_0\}.$$

\begin{lemma}
The transform $\Phi$ gives an isomorphism between the schemes $\U$ and
$\V$.
\end{lemma}

\begin{pf}
For any point $E\in \U$, $E$ is $\Psi$-WIT$_1$ and $\Hat
E\in \V$. Suppose now that $F\in \V$ and put $E=\Hat{F}$. Claim that
$E\in \U$. By Lemma \ref{whoop} the restriction of $E$ to the
general fibre of $X$ is stable, so it is only neccesary to check that
$E$ is torsion-free. Suppose $E$ has a torsion subsheaf $T$. Then since
$E$ is $\Psi$-WIT$_1$, the long exact sequence \ref{long} implies that
$T$ is $\Psi$-WIT$_1$ also, hence, by Lemma \ref{fibre}, a fibre
sheaf. Applying $\Psi$ gives a sequence
$$0\lra\Psi^0(E/T)\lRa{f} \Hat{T}\lra F\lra\Psi^1(E/T)\lra 0.$$
Since $F$ is torsion-free and $\Hat T$ is a fibre sheaf, $f$ must be an isomorphism. But,
by the result of 6.2,
$\Psi^0(E/T)$ is $\Phi$-WIT$_1$, and $\Hat{T}$ is $\Phi$-WIT$_0$. It
follows that both sheaves are zero, so $T=0$ and $E$ is torsion-free.

The proof of the lemma is completed by appealing to the general result
quoted in section 2.4.
\qed
\end{pf}

Clearly, we need to show that $\U$ and $\V$ are non-empty. Take $F\in\NN$. Then $F=L\tensor\I_Z$,
with $L\in\Pic^{\oo}(Y)$ and $Z$ a zero-dimensional subscheme of $Y$
of length $t$. By Lemma \ref{new}, $F$ is $\Phi$-WIT$_0$ precisely
when there is no non-zero map $F\to\QQ_x$ for any $x\in X$. Since
$\QQ_x$ is supported on the fibre $Y_{\pi(x)}$ of $\Hat{\pi}$, any
map $F\to\QQ_x$
factors via $F|_{Y_{\pi(x)}}$, and hence via a stable, pure dimension 1 sheaf on $Y$
of Chern class $(0,f,s)$, where $s$ is the number of points of $Z$
lying on the fibre $Y_{\pi(x)}$. Now $\QQ_x$ has Chern class
$(0,af,r)$, and is stable, so if $s<r/a$, any map $F\to\QQ_x$ is zero. This
argument, and the fact that $r>a$, gives the following results.

\begin{lemma}
Let $F=L\tensor\I_Z$, with $L\in\Pic^{\oo}(Y)$ and $Z$ a set of $t$
points $\{y_1,\cdots,y_n\}$ lying on
distinct fibres of $\pi:Y\to C$. Then $F$ is an element of $\V$.
\qed
\end{lemma}

\begin{lemma}
If $r>at$ then $\V=\NN$.
\qed
\end{lemma}

\begin{remark}
Applying $\Phi$ to
the short exact sequence
$$0\lra F\lra L\lra \OO_Z\lra 0,$$
gives a sequence
$$0\lra \Hat F\lra {\Hat L}\lra \oplus\PP_{y_i}\lra 0.$$
Let us assume for simplicity that $X$ is simply-connected. Then we see that an open
subset of $\M$ is obtained from the fixed bundle
$\Hat{\OO_Y}$ by taking $t$ distinct non-singular fibres $\{f_1,\cdots,f_t\}$ of
$\pi$ and stable bundles $P_i$ of
rank $a$ and degree 
$b$ on $f_i$, and taking the kernel of
the unique morphism
$$\Hat{\OO_Y}\lra\oplus P_i.$$
Furthermore, when $X$ is nodal, the proof of [7], Proposition III.3.11
shows that $\Hat{\OO_Y}$ is the unique sheaf (up to twists) on $X$ whose restriction
to every reduction of a fibre of $\pi$ is stable.

In the rank $2$ case, this corresponds to Friedman's method of
constructing bundles using
elementary modifications ([7], section III.3).
\end{remark}

\subsection{}
To complete the proof of Theorem 1.1 we must show that
$\M$ is
irreducible, i.e. that $\M$ has only one connected component. Let us suppose, for
contradiction, that there is a connected component $\W$ of $\M$ which
does not meet $\U$.

Let $E$ be a point of $\W$. Then $E$ is $\Psi$-WIT$_1$, and the transform $\Hat E$
is a sheaf of Chern class $(1,0,t)$ on $Y$, with a non-zero torsion
sheaf. By the argument of Lemma \ref{whoop}, the restriction of $\Hat E$ to
the general fibre of $\Hat{\pi}$ is simple, hence stable.

\begin{lemma}
Let $n\geq 1$ be an integer. Then for a general zero-dimensional subscheme $Z\in\Hilb^{rn}(Y)$, there is
a unique morphism $\Hat E\to\OO_Z$. Furthermore, for general $Z$, this
morphism surjects and the kernel $K$ is $\Phi$-WIT$_0$. The transform
$\Hat K$ is then an element of the moduli space
$$\tilde{\M}=\M_X(r,\De-(rna)f,k+rnb-rna(\De\cdot f)).$$
\end{lemma}

\begin{pf}
We may suppose that $Z$ consists of $rn$ points lying on distinct non-singular fibres
$f_1,\cdots ,f_{rn}$ of
$\Hat{\pi}$. We can also suppose that $\Hat E$ is locally-free at
each of the points of $Z$. Then there is a unique morphism $\Hat
E\to\OO_Z$ and this map surjects, giving an exact sequence
$$0\lra K\lra \Hat E\lra \OO_Z\lra 0.$$
By Lemma \ref{new}, to prove that
$K$ is $\Phi$-WIT$_0$, we must show that there are no non-zero morphisms
$K\to\QQ_x$ for any $x\in X$. We only need to check this when $\QQ_x$
is supported on one of the fibres $f_1,\cdots ,f_{rn}$ since the
restrictions of $\Hat E$ and $K$ to any other fibre are identical, and
$\Hat E$ is $\Phi$-WIT$_0$. But
we can always take $Z$ so that the restriction of $\Hat E$ to each of
the fibres $f_i$ is a degree 0 line bundle. This will be enough since
$\QQ_x$ is stable of degree $-r$.
\qed
\end{pf}
 
Twisting by
$\OO_X(anf)$ gives an isomorphism between the spaces $\tilde{\M}$ and
$\M_X(r,\De,k+n)$, so a general theorem of Gieseker-Li and O'Grady (see e.g. [16]) implies that
for large enough $n$, $\tilde{\M}$ is irreducible. It follows from
the results of 7.2 that the general element of
$\tilde{\M}$ has torsion-free transform. Now the construction of the lemma gives a rational map
$$\theta:\W\times\Hilb^{rn}(Y)\dashrightarrow\tilde{\M},$$
and since all points in the image of $\theta$ have
non-torsion-free transforms, $\theta$ cannot be dominant. But we shall
show below that the general fibre of $\theta$ is
zero-dimensional. Since $\theta$ is a map between
two varieties of the same dimension, this will give a contradiction.

Take an element of $\tilde{\M}$, and let $K$ be its
transform. We must show that there are only finitely many pairs
$$(E,Z)\in\W\times\Hilb^{rn}(Y),$$
such that $Z$ consists of $rn$ distinct points at which $\Hat E$ is
locally free, and $K=\Hat E\tensor\I_Z$.

Given such a pair, note that $Z$ does not meet the support of the
torsion subsheaf of $\Hat E$, so the torsion subsheaves $T$ of $\Hat E$
and $\Hat E\tensor\I_Z=K$ are equal. Thus $Z$ is a subset of the
finite set of points at which $K/T$ is not locally-free. This implies
that the number of possible choices of $Z$ is finite.

Finally, if we have two
pairs $(E_1,Z)$, and $(E_2,Z)$ then $E_1=E_2$, because there is only one extension
of $K$ by $\OO_Z$ which is locally-free at each of the points of $Z$.
This completes the proof.

\subsection{}
We conclude with the following simple corollary of Theorem
\ref{moduli}.

\begin{cor}
Given integers $a>0$ and $b$ with $a\m_X$ coprime to $b$, and a positive integer $t$, there are
polarizations of $X$ such that a component of the moduli space
of torsion-free stable sheaves on $X$ is isomorphic to
$$\Pic^{\oo}(\J_X(a,b))\times\Hilb^t(\J_X(a,b)).$$
\end{cor}

\begin{pf}
Take $r>at$ such that $br$ is congruent to 1 modulo $a\m_X$. Then
there exists a divisor
$\Lambda$ on $X$ such that $d=\Lambda\cdot f$ is coprime to $r$ and
$br-ad=1$. Adding multiples of
$f$ to $\Lambda$ if neccesary, it is easy to check that one can choose $k$
such that \ref{first} holds.
Applying Theorem
\ref{moduli} then gives the result.
\qed
\end{pf}

\section*{References}

\small

{EGA} {\it A. Grothendieck, J. Dieudonn{\'e},} El{\'e}ments de
g{\'e}om{\'e}trie alg{\'e}brique, Publ. Math. I.H.E.S.

[1] {\it V.I. Artamkin,} On deformations of sheaves, Math. USSR-Izv.
{\bf 32} (1989), 663-668.

[2] {\it M.F. Atiyah,} Vector bundles over an elliptic curve, Proc.
London Math. Soc. {\bf 7} (1957), 414-452.

[3] {\it W. Barth, C. Peters, A. Van de Ven,} Compact Complex
Surfaces, Ergebnisse Math. Grenzgeb. (3), vol. 4, Springer-Verlag,
1984.

[4] {\it C. Bartocci, U. Bruzzo, D. Hern{\'a}ndez-Ruip{\'e}rez,} A
Fourier-Mukai transform for stable bundles on K3 surfaces, J. reine
angew. Math. {\bf 486} (1997), 1-16.

[5] {\it A.I. Bondal,} Representations of associative algebras and
coherent sheaves, Math. USSR Izv. {\bf 34} (1990), 23-42.

[6] {\it A.I. Bondal, D.O. Orlov,} Semiorthogonal decomposition for
algebraic varieties, Preprint alg-geom 9506012.

[7] {\it R. Friedman,} Vector bundles and SO(3)-invariants for elliptic
surfaces, J. Amer. Math. Soc. {\bf 8} (1995), 29-139.
 
[8] {\it R. Friedman, J.W. Morgan,} Smooth four-manifolds and complex
surfaces, Ergebnisse Math. Grenzgeb. (3), vol. 27, Springer-Verlag,
1994.

[9] {\it R. Hartshorne,} Residues and duality, Lect. Notes Math. {\bf
20}, Springer-Verlag, 1966.

[10] {\it A. Maciocia,} Generalized Fourier-Mukai transforms, J. reine
angew. Math. {\bf 480} (1996), 197-211.

[11] {\it S. Mukai,} Duality between $D(X)$ and $D(\Hat X)$ with its
application to Picard sheaves, Nagoya Math. J. {\bf 81} (1981),
153-175.

[12] {\it S. Mukai,} Symplectic structure of the moduli space of
sheaves on an abelian or K3 surface, Invent. math. {\bf 77} (1984), 101-116.

[13] {\it S. Mukai,} On the moduli space of bundles on K3 surfaces I,
in: Vector Bundles on Algebraic Varieties, M.F. Atiyah et al., Oxford
University Press (1987), 341-413.

[14] {\it S. Mukai,} Fourier functor and its application to the moduli
of bundles on an abelian variety, Adv. Pure Math. {\bf 10} (1987),
515-550.

[15] {\it K.G. O'Grady,} The weight-two Hodge structure of moduli
spaces of sheaves on a K3 surface, Preprint alg-geom 9510001.

[16] {\it K.G. O'Grady,} Moduli of vector bundles on surfaces,
Preprint alg-geom 9609015.

[17] {\it D.O. Orlov,} Equivalences of derived categories and K3
surfaces, Preprint alg-geom 9606006.

[18] {\it C.T. Simpson,} Moduli of representations of the fundamental
group of a smooth projective variety I, Publ. Math. I.H.E.S. {\bf 79}
(1994) 47-129.

[19] {\it Y.T. Siu, G. Trautmann,} Deformations of coherent
analytic sheaves with compact supports, Memoirs of Amer. Math. Soc.
{\bf 29} (1981), no. 238.

\bigskip

Department of Mathematics and Statistics, The University of Edinburgh,
King's Buildings, Mayfield Road, Edinburgh, EH9 3JZ, UK.

email: tab@maths.ed.ac.uk

\end{document}